\documentclass{gtmon_a}
\pdfoutput=1
\usepackage{amscd,mathrsfs,float}
\usepackage[all]{xy}

\usepackage[T1,T5]{fontenc}

\proceedingstitle{Proceedings of the School and Conference in Algebraic
Topology (The Vietnam National University, Hanoi, 9-20 August 2004)}
\conferencestart{9 August 2004}
\conferenceend{20 August 2004}
\conferencename{School and Conference in Algebraic Topology}
\conferencelocation{Vietnam National University, Hanoi, Vietnam}

\editor{John Hubbuck}
\givenname{John}
\surname{Hubbuck}

\editor{Nguy\~\ecircumflex{}n H V H\uhorn{}ng}
\givenname{H\uhorn{}ng}
\surname{Nguy\~\ecircumflex{}n}

\editor{Lionel Schwartz}
\givenname{Lionel}
\surname{Schwartz}

\title{Configurations and parallelograms\\associated to centers of mass}

\author{F\,R Cohen}
\givenname{F\,R}
\surname{Cohen}
\address{Department of Mathematics\\
University of Rochester\\\newline
Rochester NY 14627\\
USA}
\email{cohf@math.rochester.edu}
\urladdr{}

\author{Y Kamiyama}
\givenname{Yasuhiko}
\surname{Kamiyama}
\address{Department of Mathematics\\
University of the Ryukyus\\\newline
Nishihara-Cho\\
Okinawa 903-0213\\
Japan}
\email{kamiyama@sci.u-ryukyu.ac.jp}
\urladdr{}

\volumenumber{11}
\issuenumber{}
\publicationyear{2007}
\papernumber{2}
\startpage{17}
\endpage{32}

\doi{}
\MR{}
\Zbl{}

\arxivreference{}


\received{6 January 2005}
\revised{15 November 2005}
\accepted{12 December 2005}
\published{14 November 2007}
\publishedonline{14 November 2007}
\proposed{}
\seconded{}
\corresponding{}
\editor{}
\version{}

\makeatletter
\def\cnewtheorem#1[#2]#3{\newtheorem{#1}{#3}[section]
\expandafter\let\csname c@#1\endcsname\c@subsection}

\let\xysavmatrix\xymatrix
\def\xymatrix{\disablesubscriptcorrection\xysavmatrix}
\AtBeginDocument{\let\bar\wbar\let\tilde\wtilde\let\hat\what}

\theoremstyle{plain}
\cnewtheorem{prop}[subsection]{Proposition}
\cnewtheorem{thm}[subsection]{Theorem}
\cnewtheorem{lem}[subsection]{Lemma}
\cnewtheorem{cor}[subsection]{Corollary}

\theoremstyle{definition}
\cnewtheorem{remarks}[subsection]{Remarks}
\cnewtheorem{remark}[subsection]{Remark}
\cnewtheorem{defn}[subsection]{Definition}
\cnewtheorem{exm}[subsection]{Example}

\makeatother  

\numberwithin{equation}{section}

\makeop{Conf}

\begin{document}

\begin{htmlabstract}
<p class="noindent">
The purpose of this article is to</p>
<p class="noindent">
1. define M(t,k) the t&ndash;fold center of mass arrangement for k
     points in the plane,<br>
2. give elementary properties of M(t,k) and<br>
3. give consequences concerning the space M(2,k) of k distinct points
in the plane, no four of which are the vertices of a parallelogram.</p>
<p class="noindent">
The main result proven in this article is that the classical
unordered configuration of k points in the plane is not a retract
up to homotopy of the space of k unordered distinct points in the
plane, no four of which are the vertices of a parallelogram. The
proof below is homotopy theoretic without an explicit computation of
the homology of these spaces.
</p>
<p class="noindent">
In addition, a second, speculative part of this article arises from
the failure of these methods in the case of odd primes p. This
failure gives rise to a candidate for the localization at odd primes
p of the double loop space of an odd sphere obtained from the
p&ndash;fold center of mass arrangement. Potential consequences are
listed.
</p>
\end{htmlabstract}

\begin{abstract}
The purpose of this article is to
{\begin{enumerate}\leftskip 20pt\rightskip25pt
\item define $M(t,k)$ the $t$--fold center of mass arrangement
  for $k$ points in the plane,
\item  give elementary properties of $M(t,k)$ and
\item give consequences concerning the space $M(2,k)$
  of $k$ distinct points in the plane, no four of which are the
  vertices of a parallelogram.
\end{enumerate}\par}
The main result proven in this article is that the classical
unordered configuration of $k$ points in the plane is not a retract
up to homotopy of the space of $k$ unordered distinct points in the
plane, no four of which are the vertices of a parallelogram. The
proof below is homotopy theoretic without an explicit computation of
the homology of these spaces.

In addition, a second, speculative part of this article arises from
the failure of these methods in the case of odd primes $p$. This
failure gives rise to a candidate for the localization at odd primes
$p$ of the double loop space of an odd sphere obtained from the
$p$--fold center of mass arrangement. Potential consequences are
listed.
\end{abstract}

\maketitle

\section{Introduction and statement of results}

Fix integers $k$ and $t$. The $t$--fold center of mass arrangement
$M(t,k)$ for integers $t$ with $ k \geq t \geq 1 $ is defined as the
subspace of the $k$--fold product $\mathbb C^k$ given by ordered
$k$--tuples of points $(x_1,\ldots,x_k)$ such that the centroids of any
set of $t$ elements in the underlying set $\{x_1,\ldots,x_k\}$
$$\sigma_t(x_{i_1},x_{i_2},\ldots,x_{i_t}) = (1/t)(
x_{i_1}+x_{i_2}+\cdots+x_{i_t})$$ are distinct for all distinct
subsets $\{x_{i_1},x_{i_2},\ldots,x_{i_t}\}$, and
$\{x_{j_1},x_{j_2},\ldots,x_{j_t}\}$. In particular, $M(t,k)$ is the
complement of the union of the hyperplanes specified by
$$\sigma_t(x_{i_1},x_{i_2},\ldots,x_{i_t})-\sigma_t(x_{j_1},x_{j_2},\ldots,x_{j_t})
= 0$$ for all pairs of unequal sets $S_I=
\{x_{i_1},x_{i_2},\ldots,x_{i_t}\}$, and $S_J=
\{x_{j_1},x_{j_2},\ldots ,x_{j_t}\}$. Write $$|S_J|$$ for the
cardinality of the set $S_J$. In case $ k < t$, define $M(t,k)$ to
be the Fadell--Neuwirth configuration space $\Conf(\mathbb C,k)$ of
ordered $k$ tuples of distinct points in $\mathbb C$ (see Fadell--Neuwirth
\cite{FN}).

Finite unions of complex hyperplanes in complex $k$--space are known
as complex hyperplane arrangements in Orlik--Terao \cite{OrlikTerao}. The space
$M(t,k)$ is a complement of a complex hyperplane arrangement.
Consider an equivalent formulation of $M(t,k)$ as the complement of
the variety $V(t,k)$ of ordered $k$--tuples $(x_1,...,x_k)$ defined
by the equation $$\prod_{S_I \neq S_J,|S_I|=|S_J|=t}
([x_{i_1}+x_{i_2}+\cdots +x_{i_t}]-[x_{j_1}+ x_{j_2}+\cdots
+x_{j_t}]) = 0$$ with $$M(t,k) = \mathbb {C}^k - V(t,k).$$

Modifications of the $M(t,k)$, $M'(t,k)$, are defined as follows:
$$M'(t,k) = \cap_{1 \leq s \leq t}M(s,k).$$ Thus $M'(t,k)$ is the
complement of the variety $W(t,k)$ of ordered $k$--tuples
$(x_1,...,x_k)$ defined by the equation
$$\prod_{S_I \neq S_J,1<q = |S_I|= |S_J|\leq t}
  \bigl([x_{i_1}+x_{i_2}+ \cdots +x_{i_q}]-[x_{j_1}+
  x_{j_2}+ \cdots +x_{j_q}]\bigr) = 0$$
with
$$M'(t,k) = \mathbb C^k - W(t,k).$$
Similarly, if $ k < t$, define $M'(t,k)$ to be $\Conf(\mathbb C,k)$.

In addition, there are natural inclusions
\[
\begin{CD}
M'(t,k) @>{}>> M(t,k) @>{}>> \Conf(\mathbb C,k).\\
\end{CD}
\]  These inclusions are equivariant with
respect to the natural action of the symmetric group on $k$ letters,
$\Sigma_k$.

Consider the $t$--fold symmetric product $\mathbb C^t/ \Sigma_t$, and
notice that there is a map $$\chi_t\co  \Conf(\mathbb C,k) \to (\mathbb
C^t/\Sigma_t)^{k \choose t}$$ gotten by choosing all $t$--element
subsets out of a set of cardinality $k$ with a fixed ordering of the
subsets. The map $\chi_t$ is given on the level of points by the
formula $$\chi_t(z_1,z_2,\ldots, z_k) = \prod_{i_1< i_2 < \ldots <
i_t} [ z_{i_1},z_{i_2},\ldots, z_{i_t}]$$ for which the points $[
z_{i_1},z_{i_2},\ldots, z_{i_t}]$ in $\mathbb C^t/\Sigma_t$ are
ordered in the product left lexicographically by indices and over
all subsets of cardinality $t$ in the set $\{z_1,z_2,\ldots, z_k\}$.

Notice that the map $\chi_t\co  \Conf(\mathbb C,k) \to (\mathbb
C^t/\Sigma_t)^{k \choose t}$ takes values in the configuration space
$ \Conf(\mathbb C^t/\Sigma_t, {k \choose t})$. Thus in what follows
below $\chi_t$ will be regarded as a map
$$\chi_t\co  \Conf(\mathbb C,k) \to \Conf\bigl(\mathbb C^t/\Sigma_t,
\tbinom{k}{t}\bigr).$$
Addition of complex numbers provides a map
$$\oplus_t \co  \mathbb C^t/\Sigma_t \to \mathbb C $$ with
$$\oplus_t ([z_1,\ldots,z_t])= z_1+ \cdots +z_t.$$
There is an induced map $$\Theta_t\co  \Conf(\mathbb C,k) \to {\mathbb
C}^{k \choose t}$$ given by the composite
\[
\begin{CD}
  \Conf( \mathbb C,k) @>{\chi_t}>> (\mathbb C^t/\Sigma_t)^ {k \choose t} @>{
  (\oplus_t)^{k \choose t} }>> {\mathbb C}^{k \choose t}.     \\
\end{CD}
\] Thus $$\Theta_t(z_1,z_2,\ldots, z_k) = \prod_{i_1< i_2 < \cdots
< i_t} ( z_{i_1}+ z_{i_2}+ \cdots + z_{i_t})$$ in ${\mathbb C}^{k
\choose t}$.

The next proposition, a useful observation, is recorded next where
$$j\co \Conf\bigl(\mathbb C, \tbinom{k}{t}\bigr) \to {\mathbb C}^{k
\choose t}$$ is
the natural inclusion. This observation is the starting point of the
results here, and provides the basic motivation for considering the
center of mass arrangement.
\begin{prop}\label{prop:proposition cartesian}
The following diagram is a pull-back (a cartesian diagram):
\[
\begin{CD}
M(t,k)       @>>>               \Conf(\mathbb C, {k \choose t})        \\
@VVV                    @VV\text{j}V                    \\
\Conf(\mathbb C,k) @>{\Theta_t}>> {\mathbb C}^{k \choose t}\\
\end{CD}
\]
\end{prop}

Notice that $M(2,k)$ is the space of ordered $k$--tuples of distinct
points such that no four of the points are the vertices of a
possibly degenerate parallelogram. Consider the natural inclusion
$M(2,k) \to \Conf(\mathbb C,k)$ modulo the action of $\Sigma_k$ the
symmetric group on $k$ letters
$$i(2,k)\co M(2,k)/\Sigma_k \to \Conf(\mathbb C,k)/\Sigma_k.$$
One question is whether there is a cross-section up to homotopy, or
even a $2$--local stable cross-section up to homotopy for this
inclusion. This last question concerns plane geometry and whether
the configuration space of distinct unordered $k$--tuples of points
in the plane can be deformed to the subspace of points, no four of
which are the vertices of a parallelogram.

\begin{thm} \label{thm:no retraction}
If $k \geq 4$, the natural map
$$i(2,k)\co M(2,k)/\Sigma_k \to \Conf(\mathbb C,k)/\Sigma_k $$ does not admit
a surjection in mod--$2$ homology, and thus does not admit a
cross-section (or a stable $2$--local cross-section) up to homotopy.
\end{thm}

The proof, homotopy theoretic without a specific computation of the
homology of these spaces, gives features of the topology of double
loop spaces which forces the maps $i(2,k)\co M(2,k)/\Sigma_k \to
\Conf(\mathbb C,k)/\Sigma_k $ for $k \geq 4$ to fail to be
epimorphisms in mod--$2$ homology.  The analogous methods applied to
the natural inclusion $$i(p,k)\co  M(p,k)/\Sigma_k \to \Conf(\mathbb
C,k)/\Sigma_k $$ for $p$ an odd prime fail to produce a non-trivial
obstruction to the existence of a stable $p$--local section. Hence a
problem unsolved here is whether $i(p,k)$ admits a stable $p$--local
cross-section. The failure of the methods here in case $p$ is an odd
prime leads to the speculation in section $2$ here concerning the
localization of the double loop space of a sphere at an odd prime
$p$.

Further properties of these arrangements are noted next. The natural
``stabilization" map for configuration spaces fails to preserve the
spaces $M(t,k)$. However, there are stabilization maps for the
modified center of mass arrangements $$S\co  M'(t,k) \to M'(t,k+1)$$
defined by $$S(x_1,\ldots,x_k) = (x_1,\ldots,x_k,\vec z)$$ where
$\vec z$ is the vector $(L,0)$ with $ L = 2t( 1 + max_{k \geq i \geq
1} ||x_i||)$. Notice that $S$ takes values in $M'(t,k+1)$, but that
the analogous map out of $ M(t,k)$ takes values in $\Conf(\mathbb
C,k)$, but not in the subspace $M(t,k+1)$.

The next result follows directly from Cohen \cite{C2} and
Cohen--May--Taylor \cite{CMT}.
\begin{thm} \label{thm: splitting}
The map $$S\co  M'(t,k) \to M'(t,k+1)$$ extends to a map
$$S_*\co  M'(t,k)\times_{\Sigma_k}Y^k \to\
M'(t,k+1)\times_{\Sigma_{k+1}}Y^{k+1}$$ which admits a stable left
inverse for any path-connected CW-complex $Y$, and thus induces a
split monomorphism in homology with any field coefficients.
\end{thm}

\begin{cor}\label{cor:sign.rep.splitting} The map $S\co  M'(t,k)/\Sigma_k \to\
M'(t,k+1)/ \Sigma_{k+1}$ induces a split mono\-morphism in homology
with coefficients in any graded permutation representation of
$\Sigma_k$, and thus by specialization to either coefficients given
by the trivial representation or the sign representation.
\end{cor}

Connections to homotopy theory as well as the motivation for
considering the spaces $M(t,k)$ and the map
$$\chi_t\co  \Conf(\mathbb C,k) \to (\mathbb C^t/\Sigma_t)^{k \choose
t}$$ defined earlier in this section are given next. These
connections arise from stable homotopy equivalences
$$H\co \Omega^2 \Sigma^2(X) \to \vee_{0 \leq
k}D_k(\Omega^2\Sigma^2(X))$$ in case $X$ is a path-connected
$CW$--complex originally proven by Snaith \cite{Snaith} and subsequently by
Cohen, May and Taylor \cite{CMT,C2} for which $D_k(\Omega^2\Sigma^2(X))$
is defined in
\fullref{speculation} here.

This stable homotopy equivalence is obtained by adding maps given by
$$h_k\co \Omega^2 \Sigma^2(X) \to \Omega^{2k}
\Sigma^{2k}D_k(\Omega^2\Sigma^2(X))$$ as observed in the appendix of
\cite{C2}. These maps do not compress through $$\Omega^{2k-1}
\Sigma^{2k-1}D_k(\Omega^2\Sigma^2(X))$$ in case $k = 2^t$, and
spaces are localized at the prime $2$ (see Cohen and Mahowald \cite{CM}).

Specialize $h_k$ to $k=p$ an odd prime and $X = S^{2n-1}$. The
spaces $M(p,k)$ and $M'(p,k)$ as well as the map $\chi_t\co 
\Conf(\mathbb C,k) \to (\mathbb C^t/\Sigma_t)^{k \choose t}$ are
introduced here in order to attempt to compress the maps
$$h_p\co \Omega^2 \Sigma^2(S^{2n-1}) \to \Omega^{2p}
\Sigma^{2p}D_p(\Omega^2\Sigma^2(S^{2n-1}))$$ through some choice of
map $$\bar h_p\co \Omega^2 \Sigma^2(S^{2n-1}) \to \Omega^{2}
\Sigma^{2}D_p(\Omega^2\Sigma^2(S^{2n-1})).$$
The map $h_p$ as given in \cite{C2,CMT} is induced on the level of
certain combinatorial models by the composite
\[
\begin{CD}
  \Conf( \mathbb C,k) @>{\chi_p}>> \Conf\bigl( \mathbb C^p/\Sigma_p ,{k
\choose p} \bigr) 
@>{\text{inclusion}}>>
  (\mathbb C^p/\Sigma_p)^ {k \choose p}.     \\
\end{CD}
\] A space $M_p(\mathbb C, X)$ together with a map
$$I_p\co M_p(\mathbb C, X) \to \Omega^2 \Sigma^2(X)$$ will be defined in
\fullref{speculation} in which configuration spaces
$\Conf(\mathbb C,k)$ used in combinatorial models of $\Omega^2
\Sigma^2(X)$ are replaced by the spaces $M(p,k)$. Furthermore, there
are continuous maps $$ h_p\co  M_p(\mathbb C, X) \to
\Omega^{2}\Sigma^{2}(D_p(\Omega^2\Sigma^2(X))).$$

It is natural to compare the homotopy types of $\Omega^2S^{2n+1}$
and $M_p(\mathbb C, S^{2n-1})$ after localization at an odd prime
$p$ by the following theorem in which $$E\co \Sigma^2(Y) \to
\Omega^{2p-2}\Sigma^{2p}(Y)$$ denotes the classical suspension map.
\begin{thm} \label{thm:comparison}
There is a commutative diagram
\[
\begin{CD}
M_p(\mathbb C, X)  @>{h_p}>> \Omega^{2}\Sigma^{2}(D_p(\Omega^2\Sigma^2(X))) \\
  @VV{I_p}V          @VV{\Omega^{2}(E)}V \\
C(\mathbb C, X) @>{h_p}>>
\Omega^{2p}\Sigma^{2p}(D_p(\Omega^2\Sigma^2(X))).
\end{CD}
\]  Thus if the map
$$I_p\co M_p(\mathbb C, S^{2n-1})  \to \Omega^2S^{2n+1}$$ is a
$p$--local equivalence, then there is a $p$--local map
$$ \bar h_p\co  \Omega^2S^{2n+1} \to
\Omega^{2}\Sigma^{2}(D_p(\Omega^2\Sigma^2(S^{2n-1})))$$ which is a
compression of the map $h_p\co \Omega^2 \Sigma^2(S^{2n-1}) \to
\Omega^{2p} \Sigma^{2p}D_p(\Omega^2\Sigma^2(S^{2n-1}))$ and which
induces an isomorphism on $H_{2np-2}(-;\mathbb F_p)$.
\end{thm}

Some consequences of the existence of $\bar h_p$ are discussed in
\fullref{speculation} here. These consequences suggest that it
would be interesting to understand the behavior of the natural map
$$I_p\co  M_p(\mathbb C, S^{2n-1}) \to C(\mathbb C, S^{2n-1})$$ on the
level of mod--$p$ homology.


The authors would like to congratulate Hu\`ynh Mui on this happy
occasion of his $60$th birthday. The work here is inspired by Mui's
mathematical work on extended power constructions as well as his
interest in configuration spaces. The authors would like to thank
Nguy\~\ecircumflex{}n~H\,V~H\uhorn{}ng as well as the other organizers of this
conference.

The first named author has been supported in part
by the NSF Grant No. DMS-0072173 and CNRS-NSF Grant No. 17149.

\section{Speculation concerning the localization of the double loop space of a 
sphere
at an odd prime $p$, and applications}\label{speculation}

The main goal of this section is to point out that if the
equivariant homology of either $M(t,k)$ or $M'(t,k)$ satisfies one
statement below, then these spaces provide a method for constructing
the localization at an odd prime $p$ of the double loop space of an
odd sphere. Some consequences are also given.

Let $R[\Sigma_k]$ denote the group ring of the symmetric group over
a commutative ring $R$ with $1$, and let $\mathcal S$ denote a left
$R[\Sigma_k]$--module. Let $X$ denote a path-connected Hausdorff
space with a free, right action of the symmetric group $\Sigma_k$.
Let $H_*(X/\Sigma_k; \mathcal S)$ denote the homology of the chain
complex $ C_*(X) \otimes_{\mathbb Z[\Sigma_k]} \mathcal S$ where $
C_*(X)$ denotes the singular chain complex of $X$.

Observe that the natural inclusion $$M(t,k) \to \Conf(\mathbb C,k)
$$ induces a homomorphism $$H_*(M(t,k)/\Sigma_k; \mathcal S) \to\
H_*(\Conf(\mathbb C,k)/\Sigma_k; \mathcal S).$$ If $t$ is equal to an
odd prime $p$, and $\mathcal S$ is the coefficient module given by
$\mathbb F_p(\pm 1)$ the field of $p$--elements with the action of
$\Sigma_k$ specified by the sign representation, then one question
is to decide whether this map induces an isomorphism in mod--$p$
homology. There is no strong evidence either way, although an
affirmative answer has interesting consequences which are described
below. The analogous question for $p = 2$ fails at once by \fullref{thm:no retraction}.

The reason for the interest in these particular homology groups is
the following observation implicit in Cohen \cite{C} as follows.

\begin{thm}
Let $\mathbb F$ denote a field. For each integer $i$ greater than
$0$, there is an isomorphism
$$ \oplus_{k \geq 0} H_{i-k(2n-1)}
(\Conf(\mathbb C,k)/\Sigma_k, \mathbb F(\pm 1)) \to\
H_i(\Omega^2S^{2n+1}; \mathbb F).$$
\end{thm}

The next corollary follows at once.

\begin{cor}
Let $\mathbb F$ denote a field, and $p$ an odd prime.

\begin{enumerate}
     \item If the natural inclusion $$M(p,k) \to \Conf(\mathbb C,k)$$
induces an isomorphism $$H_*(M(p,k)/\Sigma_k; \mathbb F_p(\pm 1))
\to H_*(\Conf(\mathbb C,k)/\Sigma_k; \mathbb F_p(\pm 1)),$$ then
there are isomorphisms
$$ \oplus_{k \geq 0} H_{i-k(2n-1)}
(M(p,k)/\Sigma_k, \mathbb F_p(\pm 1)) \to H_i(\Omega^2S^{2n+1};
\mathbb F_p).$$
     \item If the natural inclusion $$M'(p,k) \to \Conf(\mathbb R^2,k)$$
     induces an isomorphism
     $$H_*(M'(p,k)/\Sigma_k; \mathbb F_p(\pm 1)) \to H_*(\Conf(\mathbb
     C,k)/\Sigma_k; \mathbb F_p(\pm 1)),$$ then there are isomorphisms
$$ \oplus_{k \geq 0} H_{i-k(2n-1)}
(M'(p,k)/\Sigma_k, \mathbb F_p(\pm 1)) \to H_i(\Omega^2S^{2n+1};
\mathbb F_p).$$
     \end{enumerate}
\end{cor}

The spaces $M(p,k)$ and $M'(p,k)$ are used next to give analogues of
labeled configuration spaces in which the configuration space itself
is replaced by a ``center of mass construction" as given above. Let
$Y$ denote a pointed space with base-point $*$ and $W$ any
topological space. Recall the labeled configuration space
$$C(W,Y)$$ given by equivalence classes of pairs $[S,f]$ where

\begin{enumerate}
     \item $S$ is a finite subset of $W$,
     \item $f\co  S \to Y$ is a function, and
     \item $[S,f]$ is equivalent to $[S-\{p\},f|_{S-\{p\}}]$ if and
     only if $f(p) = *$.
\end{enumerate}

One theorem proven by May \cite{May} is as follows.
\begin{thm} \label{thm:May}
If $Y$ is a path-connected CW-complex, then $C(\mathbb R^n,Y)$ is
homotopy equivalent to $\Omega^n\Sigma^n(Y)$
\end{thm}

Technically, May's proof does not exhibit a map between these two
spaces. There are weak equivalences on the level of May's
construction \cite{May} $\alpha\co  C_n(Y) \to \Omega^n \Sigma^n(Y)$
and the natural evaluation map $e\co  C_n(Y) \to C(\mathbb R^n,Y)$.

Furthermore, the construction $D_k(\Omega^2\Sigma^2(X))$ is homotopy
equivalent to $$\Conf(\mathbb C,k)
\times_{\Sigma_k}X^{(k)}/\Conf(\mathbb C,k)\times_{\Sigma_k} \{*\}$$
for which $X^{(k)}$ denotes the $k$--fold smash product \cite{May}.
When localized at an odd prime $p$, $D_p(\Omega^2 S^{2n+1})$ is
homotopy equivalent to a mod--$p$ Moore space $P^{2np-1}(p)$ with a
single non-vanishing reduced homology group given by $\mathbb Z/
p\mathbb Z$ in dimension $2np-2$. This last assertion follows from
the computations in Cohen \cite{C}.

\begin{defn}
Define $$M_t(\mathbb C, Y)$$ to be the subspace of $C(\mathbb C, Y)$
given by those points for which $S$ is a subset of $M(t,k)$ with
natural inclusion denoted by $I_p\co M_t(\mathbb C, Y)\to C(\mathbb C,
Y)$ and
$$M_t'(\mathbb C, Y)$$ to be the subspace of $C(\mathbb C, Y)$ given
by those points for which $S$ is a subset of $M'(t,k)$ with natural
inclusion denoted (ambiguously) by $I_p\co M_t(\mathbb C, Y)\to
C(\mathbb C, Y)$.
\end{defn}

The next statement provides a potential method for constructing the
localization at $p$ of the double loop space of an odd sphere which
also has some useful properties.
\begin{thm} Assume that $p$ is an odd prime.

\begin{enumerate}
     \item If $M(t,k) \to \Conf(\mathbb C,k) $ induces an isomorphism
$$H_*(M(t,k)/\Sigma_k; \mathbb F_p(\pm 1)) \to H_*(\Conf(\mathbb C,k)/\Sigma_k; 
\mathbb F_p(\pm 1))$$ for $t$ an odd prime $p$, then
the natural map
$$I_p\co  M_p(\mathbb C, S^{2n-1}) \to\
\Omega^2S^{2n+1}$$ induces a mod--$p$ homology isomorphism.
     \item If $M'(t,k) \to \Conf(\mathbb C,k) $ induces an isomorphism
$$H_*(M'(t,k)/\Sigma_k; \mathbb F_p(\pm 1)) \to H_*(\Conf(\mathbb 
C,k)/\Sigma_k; \mathbb F_p(\pm 1))$$ for $t$ an odd prime $p$, then
the natural map $$I_p\co  M'_p(\mathbb C, S^{2n-1}) \to\
\Omega^2S^{2n+1}$$ induces a mod--$p$ homology isomorphism.
\end{enumerate}
\end{thm}

One consequence of this last theorem is that it implies properties
of the double suspension of $E^2\co  S^{2n-1}\to \Omega^2 S^{2n+1}$
after localization at an odd prime. In particular, the next
corollary follows directly.
\begin{cor}\label{cor:consequences}
Let $p$ denote an odd prime. If either
\begin{enumerate}
\item the natural inclusion
$M(p,k) \to \Conf(\mathbb C,k)$ induces an isomorphism
$$H_*(M(p,k)/\Sigma_k; \mathbb F_p(\pm 1)) \to H_*(\Conf(\mathbb C,k)/\Sigma_k; 
\mathbb F_p(\pm 1)),$$ or
\item the natural inclusion
$M'(p,k) \to \Conf(\mathbb R^2,k)$ induces an isomorphism
$$H_*(M'(p,k)/\Sigma_k; \mathbb F_p(\pm 1)) \to H_*(\Conf(\mathbb 
C,k)/\Sigma_k; \mathbb F_p(\pm 1)),$$
\end{enumerate}
then after localization at
$p$, the mod--$p$ Moore space $$P^{2np+1}(p)$$ is a retract of
$\Sigma^2 \Omega^2 S^{2n+1}$. In that case, the following hold:
\begin{enumerate}
\item Any map  $$\alpha\co  P^{2p+1}(p) \to S^3$$
   given by an extension of $\alpha_1\co  S^{2p} \to S^3$, which
realizes the first element of order $p$ in the homotopy groups of
the $3$--sphere induces a split epimorphism on the $p$--primary
component of homotopy groups.
   \item After localization at the prime $p$, the homotopy theoretic
   fibre of the double suspension $E^2\co  S^{2n-1} \to \Omega^2 S^{2n+1}$ is the 
fibre of
  a map $\Omega^2 S^{2np+1}\to S^{2np-1}$.
\end{enumerate}
\end{cor}

\begin{remarks}\quad
\begin{enumerate}
\item The main content of \fullref{cor:consequences} is that the
``center of mass arrangements" may provide a useful way to construct
a localization at odd primes for the double loop space of an odd
sphere. \fullref{cor:consequences} follows from \fullref{thm:comparison}
as outlined in Cohen \cite{C2}.
\item In addition, \fullref{thm:no retraction} shows that these
constructions fail to give the localization at the prime $2$ of
$\Omega^2 S^{2n+1}$.
\item It would also be interesting to see whether
there are analogous properties for the center of mass arrangement
with $\mathbb C$ replaced by $\mathbb C^n$ which may provide the
localization of $\Omega^{2n} S^{2(n+k)+1}$ at an odd prime $p$.
\end{enumerate}
\end{remarks}

\section[Sketch of \ref{prop:proposition cartesian}]{Sketch of
\fullref{prop:proposition cartesian}}
\label{proof of cartesian}

Notice that the set theoretic pull-back in the diagram given in
\fullref{prop:proposition cartesian} is precisely the
subspace of the configuration space given by the $t$--fold center of
mass arrangement $M(t,k)$.

\section[Calculations at the prime 2, and the proof of \ref{thm:no
retraction}]{Calculations at the prime 2, and the proof of \fullref{thm:no
retraction}}
\label{ the prime two }

The method here of comparing the homology of the center of mass
arrangement with that of the configuration space uses some
additional topology. Here, consider the natural inclusion $M(p,k)
\to \Conf(\mathbb C,k)$ together with the induced map $$I_p\co 
M_p(\mathbb C, S^{2n-1}) \to C(\mathbb C, S^{2n-1}).$$

The space $C(\mathbb C, S^{2n-1})$ is homotopy equivalent to
$\Omega^2 S^{2n+1}$ (see May \cite{May}). In addition, the spaces $M_p(\mathbb
C, S^{2n-1})$, and $C(\mathbb C, S^{2n-1})$ admit stable
decompositions which are compatible  by the remarks in
Cohen \cite{C2} and Cohen--May--Taylor \cite{CMT}. Notice that the
inclusion $M(t,k)\to \Conf(\mathbb
C,k)$ is the identity in case $k \leq t$ by definition of $M(t,k)$.
Thus the induced maps on stable summands
$$D_j(M_p(\mathbb C, S^{2n-1})) \to D_j(C(\mathbb C, S^{2n-1}))$$
is the identity in case $j \leq p,$ a feature which is used below.

Let $p = 2$, and consider the second stable summand
$$D_2(M_2(\mathbb C, S^{2n-1})) = D_2(C(\mathbb C, S^{2n-1})).$$
This stable summand is homotopy equivalent to
$$(S^1\times_{\Sigma_2}S^{4n-2})/ (S^1 \times _{\Sigma_2}*)$$ which is itself
homotopy equivalent to  $$\Sigma^{4n-3}(\mathbb R \mathbb P^2).$$
Let $u$ denote a basis element for $H_{4n-2}(\Sigma^{4n-3}(\mathbb
R\mathbb P^2);\mathbb F_2)$, and $v$ denote a basis element for
$H_{4n-1}(\Sigma^{4n-3}(\mathbb R\mathbb P^2);\mathbb F_2)$.

In addition, there is a strictly commutative diagram
\[
\begin{CD}
M_2(\mathbb C, S^{2n-1})  @>{h_2}>> 
\Omega^{\infty}\Sigma^{\infty}(D_2(C(\mathbb C, S^{2n-1}))) \\
  @VV{I_2}V          @VV{1}V \\
C(\mathbb R^2, S^{2n-1}) @>{h_2}>>
\Omega^{\infty}\Sigma^{\infty}(D_2(C(\mathbb C, S^{2n-1})))
\end{CD}
\] which gives the fact that the space $D_2(C(\mathbb C, S^{2n-1})) = 
\Sigma^{4n-3}(\mathbb R\mathbb P^2)$
is a stable retract of both spaces, in a way which is compatible
with the natural stable decompositions.

Further, by \fullref{prop:proposition cartesian}, together
with the definition \cite{C2} of the map $$h_2\co M_2(\mathbb C,
S^{2n-1}) \to \Omega^{\infty}\Sigma^{\infty}(D_2(C(\mathbb C,
S^{2n-1}))),$$ there is a compression of this map through
$\Omega^{2}\Sigma^{2}(D_2(C(\mathbb C, S^{2n-1})))$. Thus, there is
a commutative diagram given as follows.
\[
\begin{CD}
M_2(\mathbb C, S^{2n-1})  @>{h_2}>> \Omega^{2}\Sigma^{2}(D_2(C(\mathbb C, 
S^{2n-1}))) \\
  @VV{I_2}V          @VV{\Omega^2(E)}V \\
C(\mathbb R^2, S^{2n-1}) @>{h_2}>>
\Omega^{\infty}\Sigma^{\infty}(D_2(C(\mathbb C, S^{2n-1}))).
\end{CD}
\] These remarks have the following consequence for which $Q_i(x)$ is the 
standard
notation for Araki--Kudo--Dyer--Lashof operations as described in
\cite{C}.

\begin{lem} \label{lem: compressions}
The image of the map
$$h_2\co  M_2(\mathbb C, S^{2n-1}) \to
  \Omega^{\infty}\Sigma^{\infty}(D_2(M_2(\mathbb C, S^{2n-1})))$$
on
the level of mod--$2$ homology is contained in the subalgebra
generated by the elements $x$, and $Q_1^q(x)$ for $q \geq 1$ for
which $x$ is an element of a basis for the mod--$2$ homology of
$\Sigma^{4n-3}(\mathbb R\mathbb P^2)$ given by $\{u,v\}$. In
particular, the element $Q_3(x)$ cannot appear as a non-trivial
summand of the image.
\end{lem}

\begin{lem} \label{lem: surjection}
If $k \geq 4$, and the natural map
$$M(2,k)/\Sigma_k \to \Conf(\mathbb C,k)/\Sigma_k $$
induces a surjection in mod--$2$ homology, then the natural map
$$M(2,4)/\Sigma_4 \to \Conf(\mathbb C,4)/\Sigma_4 $$ induces a
surjection in mod--$2$ homology.
\end{lem}

\begin{proof}
Notice that if $k \geq 4$,  the space $\Conf(\mathbb C,4)/\Sigma_4$
is a stable retract of the space $\Conf(\mathbb C,k)/\Sigma_k$ via a
map induced by the transfer obtained from the natural
$\Sigma_k$--cover (see Cohen--May--Taylor \cite{CMT2}). Thus there is a commutative diagram
\[
\begin{CD}
  \Sigma^{2k}(M(2,k)/\Sigma_k) @>{}>>  \Sigma^{2k}(\Conf(\mathbb C,k)/\Sigma_k) 
\\
  @VV{tr}V          @VV{tr}V \\
\Sigma^{2k}(M(2,4)/\Sigma_4) @>{}>> \Sigma^{2k}(\Conf(\mathbb
C,4)/\Sigma_4)
\end{CD}
\] in which the vertical maps are induced by the natural
transfer. Hence the natural\break map $M(2,4)/\Sigma_4 \to \Conf(\mathbb
C,4)/\Sigma_4$ induces a surjection on mod--$2$ homology\break as the maps
$M(2,k)/\Sigma_k \to \Conf(\mathbb C,k)/\Sigma_k$ as well as
$tr\co \Sigma^{2k}(\Conf(\mathbb C,k)/\Sigma_k) \to$\break$
\Sigma^{2k}(\Conf(\mathbb C,4)/\Sigma_4)$ induce surjections on
mod--$2$ homology by \cite{CMT2}.
\end{proof}

The proof of \fullref{thm:no retraction} is given next.
\begin{proof}
Assume that the natural inclusion $M(2,k) \to \Conf(\mathbb C,k)$
induces an epimorphism on the level of $H_*(M(2,k)/\Sigma_k; \mathbb
F_2) \to H_*(\Conf(\mathbb C,k)/\Sigma_k;\mathbb  F_2)$ for some $ k
\geq 4$. Then by \fullref{lem: surjection}, the natural map
$M(2,4)/\Sigma_4 \to \Conf(\mathbb C,4)/\Sigma_4$ induces a
surjection on mod--$2$ homology, and the induced map $H_*(M_2(\mathbb
C, S^{2n-1});\mathbb F_2) \to H_*(C(\mathbb C, S^{2n-1});\mathbb
F_2)$ is an epimorphism in dimensions $\leq 8n-1$. This will lead to
a contradiction.

Notice that
\begin{enumerate}
\item ${h_2}_*({x_{2n-1}}^2) = u$,

\item ${h_2}_*(Q_1(x_{2n-1})) = v$ and

\item ${h_2}_*(Q_1Q_1(x_{2n-1})) = AQ_1(v) + BQ_3(u)$ for scalars
$A$, and $B$ where $u$ is the unique non-zero class in
$H_{4n-2}(D_2(C(\mathbb R^2, S^{2n-1})); \mathbb F_2)$, and $v$ is
the unique non-zero class in $H_{4n-1}(D_2(C(\mathbb C, S^{2n-1}));
\mathbb F_2)$ (see Cohen \cite{C}).
\end{enumerate}
A direct computation using $Sq^1_*$, $Sq^2_*$, and the coproduct
gives $$A = B  =1.$$ The details are as follows. Notice that
$Sq^2_*(Q_1Q_1(x_{2n-1})) = 0$, but that
$$Sq^2_*(Q_1(v)) = Q_1(u) = Sq^2_*(Q_3(u)).$$
Thus $A = B$. Furthermore $Sq^1_*(Q_1Q_1(x_{2n-1}))
= Q_1(x_{2n-1})^2.$

Finally notice that ${h_2}_*(({x_{2n-1}}^2)\cdot Q_1(x_{2n-1})) = u
\cdot v + P$ where $P$ is a primitive element. The only non-zero
choice for this primitive element $P$ is $Q_1(u)$. However,
$Sq^1_*(P) = 0$. Thus ${h_2}_*({x_{2n-1}}^4) = Sq^1_*(u \cdot v + P)
= u^2$. Hence
$$Sq^2_*Sq^1_*{h_2}_*(Q_1Q_1(x_{2n-1}))= u^2$$ and $A= B = 1$.

It follows that if the natural map $M(2,4)/\Sigma_4 \to \Conf(\mathbb
C,4)/\Sigma_4$ induces a surjection on mod--$2$ homology, then the
class $Q_1(v) + Q_3(u)$ is in the image of the composite of the
following two maps:
$${I_2}_*\co  H_*(M_2(\mathbb C, S^{2n-1});\mathbb F_2) \to H_*(C(\mathbb C,
S^{2n-1});\mathbb F_2),$$ and
$${h_2}_*\co  H_*(C(\mathbb C, S^{2n-1});\mathbb F_2) \to
H_*(\Omega^{\infty}\Sigma^{\infty}D_2(C(\mathbb C,
S^{2n-1}));\mathbb F_2).$$
Thus the above computation gives that the class $Q_1(v) + Q_3(u)$ is
in the image of the composite
\[
\begin{CD}
H_*(M_2(\mathbb C, S^{2n-1});\mathbb F_2)  @>{{h_2}_* \circ
{I_2}_*}>> H_*(\Omega^{\infty}\Sigma^{\infty}D_2(C(\mathbb C,
S^{2n-1}));\mathbb F_2).
\end{CD}
\]
By \fullref{lem: compressions}, the class $Q_1(v) + Q_3(u)$ cannot
be in the image of the map
$$H_*(\Omega^{2}\Sigma^{2}(D_2(M_2(\mathbb C, S^{2n-1})));\mathbb F_2)
\to H_*(\Omega^{\infty}\Sigma^{\infty}(D_2(C(\mathbb C,
S^{2n-1})));\mathbb F_2).$$ Hence, (1) the natural map
$M(2,4)/\Sigma_4 \to \Conf(\mathbb C,4)/\Sigma_4$ cannot induce a
surjection on mod--$2$ homology and (2) the natural map
$M(2,k)/\Sigma_k \to \Conf(\mathbb C,k)/\Sigma_k$ for $ k \geq 4$
cannot induce a surjection on mod--$2$ homology. The theorem
follows.
\end{proof}

\section[Sketch of \ref{thm: splitting} and
\ref{cor:sign.rep.splitting}]{Sketch of \fullref{thm: splitting} and
\fullref{cor:sign.rep.splitting}}
\label{proof of splitting}
The proof of follows \fullref{thm: splitting} at once from the
constructions in the appendix of Cohen \cite{C2} or the main theorem in
Cohen--May--Taylor \cite{CMT} where it was shown that these maps admit stable right
homotopy inverses.

To check \fullref{cor:sign.rep.splitting}, notice that the
sign representation is given by the action of the symmetric group on
the top non-vanishing homology group of $(S^1)^n$. The corollary
follows from \fullref{thm: splitting}.

\section[Sketch of \ref{thm:comparison}]{Sketch of
\fullref{thm:comparison}} \label{proof of comparison}
The commutativity of the diagram in \fullref{thm:comparison} follows by
definition. That the map $$h_p\co  \Omega^2S^{2n+1} \to
\Omega^{2p}\Sigma^{2p}(D_p(\Omega^2\Sigma^2(S^{2n-1})))$$ induces an
isomorphism on $H_{2np-2}(-;\mathbb F_p)$ is checked in Cohen \cite{C2}.
Since $h_p$ induces an isomorphism on $H_{2np-2}(-;\mathbb F_p)$, it
follows from the known homology of these spaces that $\bar h_p$ does
also. Given a map with the homological properties of $\bar h_p$, the
proof of \fullref{thm:comparison} follows from \cite{C2}.

\begin{remark}
The goal of this approach is to try to desuspend a map
analogous to that given by Selick \cite{selick}.
\end{remark}

\bibliographystyle{gtart}
\bibliography{link}

\begin{thebibliography}{}
\providecommand\bibmarginpar{\leavevmode\marginpar}
\def\urlstyle#1{{\tt #1}}

\bibitem{C}
\textbf{F\,R Cohen}, \emph{The homology of $\mathcal{C}_{n+1}$--spaces}, from:
  ``The homology of iterated loop spaces'', Lecture Notes in Mathematics 533,
  Springer, Berlin (1976)  207--351

\bibitem{C2}
\textbf{F\,R Cohen}, \href{http://dx.doi.org/10.1007/BF01215483} {\emph{The
  unstable decomposition of {$\Omega \sp{2}\Sigma \sp{2}X$} and its
  applications}}, Math. Z. 182 (1983) 553--568 \xox{MR}{701370}

\bibitem{CM}
\textbf{F\,R Cohen}, \textbf{M\,E Mahowald}, \emph{Unstable properties of
  {$\Omega \sp{n}S\sp{n+k}$}}, from: ``Symposium on Algebraic Topology in honor
  of Jos\'e Adem (Oaxtepec, 1981)'', Contemp. Math. 12, Amer. Math. Soc.,
  Providence, R.I. (1982)  81--90 \xox{MR}{676319}

\bibitem{CMT}
\textbf{F\,R Cohen}, \textbf{J\,P May}, \textbf{L\,R Taylor}, \emph{Splitting
  of certain spaces {$CX$}}, Math. Proc. Cambridge Philos. Soc. 84 (1978)
  465--496 \xox{MR}{503007}

\bibitem{CMT2}
\textbf{F\,R Cohen}, \textbf{J\,P May}, \textbf{L\,R Taylor}, \emph{Splitting
  of some more spaces}, Math. Proc. Cambridge Philos. Soc. 86 (1979) 227--236
  \xox{MR}{538744}

\bibitem{FN}
\textbf{E Fadell}, \textbf{L Neuwirth}, \emph{Configuration spaces}, Math.
  Scand. 10 (1962) 111--118 \xox{MR}{0141126}

\bibitem{May}
\textbf{J\,P May}, \emph{The geometry of iterated loop spaces}, Lecture Notes
  in Mathematics 271, Springer, Berlin (1972) \xox{MR}{0420610}

\bibitem{OrlikTerao}
\textbf{P Orlik}, \textbf{H Terao}, \emph{Arrangements of hyperplanes},
  Grundlehren der Mathematischen Wissenschaften 300, Springer, Berlin (1992)
  \xox{MR}{1217488}

\bibitem{selick}
\textbf{P Selick}, \href{http://dx.doi.org/10.1016/0040-9383(78)90007-1}
  {\emph{Odd primary torsion in $\pi_k(S^3)$}}, Topology 17 (1978) 407--412
  \xox{MR}{516219}

\bibitem{Snaith}
\textbf{V\,P Snaith}, \emph{A stable decomposition of {$\Omega
  \sp{n}S\sp{n}X$}}, J. London Math. Soc. $(2)$ 7 (1974) 577--583
  \xox{MR}{0339155}

\end{thebibliography}

\end{document}